\documentclass{article}
\usepackage{amsfonts}
\usepackage{amsmath}

\newtheorem{theorem}{Theorem}

\newtheorem{corollary}[theorem]{Corollary}

\newtheorem{lemma}[theorem]{Lemma}

\newtheorem{proposition}[theorem]{Proposition}
\newtheorem{remark}[theorem]{Remark}

\newenvironment{proof}[1][Proof]{\noindent\textbf{#1.} }{\ \rule{0.5em}{0.5em}}
\input{tcilatex}

\begin{document}

\title{ON THE MEHLIG-WILKINSON\ REPRESENTATION\ OF\ METAPLECTIC\ OPERATORS}
\author{Maurice A. de Gosson}
\maketitle

\begin{abstract}
We study the Weyl representation of metaplectic operators suggested by
earlier work of Mehlig and Wilkinson. We give precise calculations for the
associated Maslov indices; these intervene in a crucial way in the
Gutzwiller formula of semiclassical mechanics.
\end{abstract}

\section{Introduction}

In an interesting paper \cite{MW} the physicists Mehlig and Wilkinson
introduce, in connection with their study of the Gutzwiller semiclassical
trace formula, a class of unitary operators $\widehat{S}:L^{2}(\mathbb{R}%
^{n})\longrightarrow L^{2}(\mathbb{R}^{n})$. These operators are defined as
follows: let $S\in Sp(n)$ have no eigenvalue equal to one; to $S$ one
associates the Weyl operator 
\begin{equation}
\widehat{R}(S)=\left( \frac{1}{2\pi }\right) ^{n}\frac{i^{\nu }}{\sqrt{|\det
(S-I)}}\int e^{\frac{i}{2}\left\langle M_{S}z_{0},z_{0}\right\rangle }%
\widehat{T}(z_{0})d^{2n}z_{0}  \label{sf3}
\end{equation}%
where $\widehat{T}(z_{0})$ is the Weyl--Heisenberg operator and%
\begin{equation}
M_{S}=\tfrac{1}{2}J(S+I)(S-I)^{-1}  \label{ms}
\end{equation}%
$I$ being the identity and $J$ the standard symplectic matrix (see below).
The index $\nu $ is an integer related to the sign of $\det (S-I)$, and
which is not studied in the general case in \cite{MW}. Mehlig and Wilkinson
moreover show that 
\begin{equation}
\widehat{R}(SS^{\prime })=\pm \widehat{R}(S)\widehat{R}(S^{\prime })
\label{sf4}
\end{equation}%
for all $S,S^{\prime }$ for which both sides are defined. They claim that
these operators belong to the metaplectic group. This property is however
not quite obvious; what is acceptably \textquotedblleft
obvious\textquotedblright\ is that $\widehat{R}(S)$ is a \textit{multiple}
by a scalar factor of modulus one of any of the two metaplectic operators $%
\pm \widehat{S}$ associated to $Mp(n)$; this is achieved using the
metaplectic covariance of the Heisenberg--Weyl operators (see below). The
purpose of this paper is to precise Mehlig and Wilkinson's statement by
comparing explicitly the integer $\nu $ in (\ref{sf3}) with the Maslov
indices on the metaplectic group we have studied in a previous work \cite%
{AIF}. This is indeed important --and not just an academic exercise-- since
the ultimate goal in \cite{MW} is to apply formula (\ref{sf3}) to give a new
proof of Gutzwiller's trace formula for chaotic systems. It is well-known
that the calculation of the associated \textquotedblleft Maslov
indices\textquotedblright\ is notoriously difficult: it suffices to have a
look on the impressive bibliography devoted to that embarrassingly subtle
topic. We will, in addition, give a semiclassical interpretation of $%
\widehat{R}(S)$, expressed in terms of the phase space wavefunctions we
introduced in \cite{Bullsci,IHP}.

\begin{remark}
An alternative approach to the results of this paper would be to use Howe's
beautiful \textquotedblleft oscillator group\textquotedblright\ method \cite%
{Howe} (see \cite{Folland} for a review); this would however in our case
lead to unnecessary technical complications.
\end{remark}

\subsection*{Notations}

We denote by $\sigma $ the canonical symplectic form on $\mathbb{R}_{z}^{2n}=%
\mathbb{R}_{x}^{n}\times \mathbb{R}_{p}^{n}$%
\begin{equation*}
\sigma (z,z^{\prime })=\left\langle p,x^{\prime }\right\rangle -\left\langle
p^{\prime },x\right\rangle \text{ \ if \ }z=(x,p)\text{, }z^{\prime
}=(x^{\prime }p^{\prime }) 
\end{equation*}%
that is%
\begin{equation*}
\sigma (z,z^{\prime })=\left\langle Jz,z^{\prime }\right\rangle \text{ \ \ ,
\ }J=%
\begin{bmatrix}
0 & I \\ 
-I & 0%
\end{bmatrix}%
\text{.} 
\end{equation*}%
The real symplectic group $Sp(n)$ consists of all linear automorphisms $S:%
\mathbb{R}_{z}^{2n}\longrightarrow \mathbb{R}_{z}^{2n}$ such that $\sigma
(Sz,Sz^{\prime })=\sigma (z,z^{\prime })$ for all $z,z^{\prime }$. It is a
connected Lie group. We denote by $\ell _{X}$ and $\ell _{P}$ the Lagrangian
planes $\mathbb{R}_{x}^{n}\times 0$ and $0\times \mathbb{R}_{p}^{n}$,
respectively. $\mathcal{S}(\mathbb{R}^{n})$ is the Schwartz space of rapidly
decreasing functions on $\mathbb{R}^{n}$, and its dual $\mathcal{S}^{\prime
}(\mathbb{R}^{n})$ the space of tempered distributions.

\section{Prerequisites}

\subsection{Standard theory of $Mp(n)$: Review}

The material of this first subsection is quite classical; see for instance 
\cite{Folland,AIF} and the references therein.

Every $S\in Mp(n)$ is the product of two \textquotedblleft quadratic Fourier
transforms\textquotedblright , which are operators $S_{W,m}$ defined on $%
\mathcal{S}(X)$ by%
\begin{equation}
S_{W,m}f(x)=\left( \frac{1}{2\pi i}\right) ^{n}i^{m}\sqrt{|\det L|}\int
e^{iW(x,x^{\prime })}f(x^{\prime })d^{n}x^{\prime }  \label{swm1}
\end{equation}%
where $W$ is a quadratic form in the variables $x,x^{\prime }$ of the type%
\begin{equation}
W(x,x^{\prime })=\frac{1}{2}\langle Px,x\rangle -\langle Lx,x^{\prime
}\rangle +\frac{1}{2}\langle Qx^{\prime },x^{\prime }\rangle  \label{wplq}
\end{equation}%
with $P=P^{T}$, $Q=Q^{T}$, $\det L\neq 0$. The integer $m$ appearing in (\ref%
{swm1}) corresponds to a choice of $\arg \det L$:%
\begin{equation*}
m\pi \equiv \arg \det L\text{ \ }\func{mod}2\pi 
\end{equation*}%
and to every $W$ there thus corresponds two different choices of $m$ modulo $%
4$: if $m$ is one choice, then $m+2$ is the other (this of course reflects
the fact that $Mp(n)$ is a two-fold covering of $Sp(n)$). The projection $%
\pi :Mp(n)\longrightarrow Sp(n)$ is entirely specified by the datum of each $%
\pi (S_{W,m})$, and we have $\pi (S_{W,m})=S_{W}$ where%
\begin{equation*}
(x,p)=S_{W}(x^{\prime },p^{\prime })\Longleftrightarrow p=\partial
_{x}W(x,x^{\prime })\text{ \ and }p^{\prime }=-\partial _{x^{\prime
}}W(x,x^{\prime })\text{.} 
\end{equation*}%
In particular, 
\begin{equation}
S_{W}=%
\begin{bmatrix}
L^{-1}Q & L^{-1} \\ 
PL^{-1}Q-L^{T} & PL^{-1}%
\end{bmatrix}
\label{plq}
\end{equation}%
is the free symplectic automorphism generated by the quadratic form $W$;
observe that $S_{W}\ell _{P}\cap \ell _{P}=0$ for every $W$. The inverse $%
\widehat{S}_{W,m}^{-1}=\widehat{S}_{W,m}^{\ast }$ of $\widehat{S}_{W,m}$ is
the operator $S_{W^{\ast },m^{\ast }}$ where $W^{\ast }(x,x^{\prime
})=-W(x^{\prime },x)$ and $m^{\ast }=n-m$, $\func{mod}4$. Note that if
conversely $S$ is a free symplectic matrix 
\begin{equation}
S=%
\begin{bmatrix}
A & B \\ 
C & D%
\end{bmatrix}%
\in Sp(n)\text{ \ , \ }\det B\neq 0  \label{free}
\end{equation}%
then $S=S_{W}$ with $P=B^{-1}A$, $L=B^{-1}$, $Q=DB^{-1}$.

\subsection{Heisenberg--Weyl operators}

For $z_{0}=(x_{0},p_{0})$ we denote by $T(z_{0})$ the translation $%
z\longmapsto z+z_{0}$; it acts on functions by push-forward: $%
T(z_{0})f(z)=f(z-z_{0})$. We denote by $\widehat{T}(z_{0})$ the
corresponding Heisenberg--Weyl operator: for $f\in \mathcal{S}(\mathbb{R}%
^{n})$ we have 
\begin{equation*}
\widehat{T}(z_{0})=e^{(\left\langle p_{0},x\right\rangle -\tfrac{1}{2}%
\left\langle p_{0},x_{0}\right\rangle )}f(x-x_{0})\text{.} 
\end{equation*}%
The operators $\widehat{T}(z_{0})$ satisfy the metaplectic covariance
formula:%
\begin{equation}
\widehat{S}\widehat{T}(z)=\widehat{T}(Sz)\widehat{S}\text{ \ \ }(S=\pi (%
\widehat{S}))  \label{meco}
\end{equation}%
for every $\widehat{S}\in Mp(n)$ and $z$. In fact, the metaplectic operators
are the only unitary operators, up to a an factor in $S^{1}$ satisfying (\ref%
{meco}):

\begin{quote}
\emph{For every }$S\in Sp(n)$ \emph{there exists a unitary transformation }$%
\widehat{U}$ in $L^{2}(\mathbb{R}^{n})$ \emph{satisfying (\ref{meco}) and }$%
\widehat{U}$ \emph{is uniquely determined apart from a constant factor of
modulus one.}
\end{quote}

The Heisenberg--Weyl operators moreover satisfy the relations

\begin{equation}
\widehat{T}(z_{0})\widehat{T}(z_{1})=e^{-i\sigma (z_{0},z_{1})}\widehat{T}%
(z_{1})\widehat{T}(z_{0})  \label{noco1}
\end{equation}

\begin{equation}
\widehat{T}(z_{0}+z_{1})=e^{-\tfrac{i}{2}\sigma (z_{0},z_{1})}\widehat{T}%
(z_{0})\widehat{T}(z_{1})  \label{noco2}
\end{equation}
as is easily seen from the definition of these operators.

\subsection{Weyl operators}

Let $a^{w}$ be the Weyl operator with symbol $a$: 
\begin{equation*}
a^{w}f(x)=\left( \tfrac{1}{2\pi }\right) ^{n}\int e^{i\left\langle
p,x-y\right\rangle }a(\tfrac{1}{2}(x+y),p)f(y)d^{n}yd^{n}p\text{; } 
\end{equation*}%
where $f\in \mathcal{S}(\mathbb{R}^{n})$equivalently%
\begin{equation*}
a^{w}=\int a_{\sigma }(z_{0})\widehat{T}(z_{0})d^{n}z_{0} 
\end{equation*}%
where $a_{\sigma }$ is the symplectic Fourier transform $F_{\sigma }a$
defined by%
\begin{equation*}
F_{\sigma }a(z)=\left( \tfrac{1}{2\pi }\right) ^{n}\int e^{i\sigma
(z,z^{\prime })}a(z^{\prime })d^{2n}z^{\prime }\text{.} 
\end{equation*}%
The kernel of $a^{w}$ is related to $a$ by the formula 
\begin{equation*}
a(x,p)=\int e^{-i\left\langle p,y\right\rangle }K(x+\tfrac{1}{2}y,x-\tfrac{1%
}{2}y)d^{n}y\text{.} 
\end{equation*}

The Mehlig--Wilkinson operator (\ref{sf3}) is the Weyl operator with twisted
Weyl symbol%
\begin{equation}
a_{\sigma }(z)=\left( \frac{1}{2\pi }\right) ^{n}\frac{i^{\nu }}{\sqrt{|\det
(S-I)}}e^{\frac{i}{2}\left\langle M_{S}z_{0},z_{0}\right\rangle }\text{.}
\label{asig}
\end{equation}

\subsection{Generalized Fresnel Formula}

We will use the following formula, generalizing the usual Fresnel integral
to complex Gaussians. Let $M$ be a real symmetric $n\times n$ matrix. If $M$
is invertible then the Fourier transform of the exponential $\exp
(i\left\langle Mx,x\right\rangle /2)$ is given by the formula%
\begin{equation}
\left( \tfrac{1}{2\pi }\right) ^{n/2}\int e^{-i\left\langle p,x\right\rangle
}e^{\frac{i}{2}\left\langle Mx,x\right\rangle }d^{n}x=|\det M|^{-1/2}e^{%
\frac{i\pi }{4}\limfunc{sgn}M}e^{-\frac{i}{2}\left\langle
M^{-1}x,x\right\rangle }  \label{fres}
\end{equation}%
where $\limfunc{sgn}M$, the \textquotedblleft signature\textquotedblright\
of $M$, is the number of $>0$ eigenvalues of $M$ minus the number of $<0$
eigenvalues.

For a proof see for instance \cite{Folland}, App. A.

\section{Discussion of the Mehlig--Wilkinson Formula}

The Mehlig--Wilkinson operators $\widehat{R}(S)$ are Weyl operators with
twisted Weyl symbol%
\begin{equation*}
a_{\sigma }(z)=\left( \frac{1}{2\pi }\right) ^{n}\frac{i^{\nu }}{\sqrt{|\det
(S-I)}}e^{\frac{i}{2}\left\langle M_{S}z_{0},z_{0}\right\rangle }\text{.} 
\end{equation*}%
We begin by giving two straightforward alternative formulations of these
operators.

\subsection{Equivalent formulations}

We begin by remarking that the matrix $M_{S}=\frac{1}{2}J(S+I)(S-I)^{-1}$ is
symmetric; this immediately follows from the conditions%
\begin{equation*}
S\in Sp(n)\Longleftrightarrow S^{T}JS=J\Longleftrightarrow SJS^{T}=J\text{.} 
\end{equation*}%
Notice that (\ref{ms}) can be \textquotedblleft solved\textquotedblright\ in 
$S$, yielding $S=(2M-J)^{-1}(2M+J)$.

\begin{proposition}
The operator 
\begin{equation}
\widehat{R}(S)=\left( \frac{1}{2\pi }\right) ^{n}\frac{i^{\nu }}{\sqrt{|\det
(S-I)}}\int e^{\frac{i}{2}\left\langle M_{S}z_{0},z_{0}\right\rangle }%
\widehat{T}(z_{0})d^{2n}z_{0}  \label{alf0}
\end{equation}%
can be written in the following alternative two forms:%
\begin{equation}
\widehat{R}(S)=\left( \frac{1}{2\pi }\right) ^{n}\frac{i^{\nu }}{\sqrt{|\det
(S-I)}}\int e^{-\frac{i}{2}\sigma (Sz_{0},z_{0})}\widehat{T}%
((S-I)z_{0})d^{2n}z_{0}  \label{alf1}
\end{equation}%
\begin{equation}
\widehat{R}(S)=\left( \frac{1}{2\pi }\right) ^{n}i^{\nu }\sqrt{|\det (S-I)}%
\int \widehat{T}(Sz_{0})\widehat{T}(-z_{0})d^{2n}z_{0}  \label{alf2}
\end{equation}
for $\det (S-I)\neq 0$.
\end{proposition}

\begin{proof}
We have 
\begin{equation*}
\tfrac{1}{2}J(S+I)(S-I)^{-1}=\tfrac{1}{2}J+J(S-I)^{-1} 
\end{equation*}%
hence, in view of the antisymmetry of $J$,%
\begin{equation*}
\left\langle M_{S}z_{0},z_{0}\right\rangle =\left\langle
J(S-I)^{-1}z_{0},z_{0}\right\rangle =\sigma ((S-I)^{-1}z_{0},z_{0}) 
\end{equation*}%
Performing the change of variables $z_{0}\longmapsto (S-I)^{-1}z_{0}$ we can
rewrite the integral in the right hand side of (\ref{alf0}) as%
\begin{eqnarray*}
\int e^{\frac{i}{2}\left\langle M_{S}z_{0},z_{0}\right\rangle }\widehat{T}%
(z)d^{2n}z_{0} &=&\int e^{\frac{i}{2}\sigma (z_{0},(S-I)z_{0})}\widehat{T}%
((S-I)z_{0})d^{2n}z_{0} \\
&=&\int e^{-\frac{i}{2}\sigma (Sz_{0},z_{0})}\widehat{T}%
((S-I)z_{0})d^{2n}z_{0}
\end{eqnarray*}%
hence (\ref{alf1}). Taking into account the relation (\ref{noco2}) we have%
\begin{equation*}
\widehat{T}((S-I)z_{0})=e^{-\tfrac{i}{2}\sigma (Sz_{0},z_{0})}\widehat{T}%
(Sz_{0})\widehat{T}(-z_{0}) 
\end{equation*}%
and formula (\ref{alf2}) follows.
\end{proof}

\begin{corollary}
We have $\widehat{R}(S)=c_{S}\widehat{S}_{W,m}$ where $c$ is a complex
constant with $|c|=1$.
\end{corollary}

\begin{proof}
We begin by noting that $\widehat{R}(S)$ satisfies the metaplectic
covariance relation 
\begin{equation*}
\widehat{R}(S)\widehat{T}(z_{0})=\widehat{T}(Sz_{0})\widehat{R}(S) 
\end{equation*}%
as immediately follows from the alternative form (\ref{alf2}) of $\widehat{R}%
(S)$. On the other hand, a straightforward calculation using formula (\ref%
{alf1}) shows that $\widehat{R}(S)$ is unitary, hence the claim.
\end{proof}

\subsection{The case $\widehat{S}=\widehat{S}_{W,m}$}

We are going to show that the Mehlig--Wilkinson operators coincide with the
metaplectic operators $\widehat{S}_{W,m}$ when $S=S_{W}$ and we will
thereafter determine the correct choice for $\nu $; we will see thast it is
related by a simple formula to the usual Maslov index as defined in \cite%
{AIF}.

Let us first prove the following technical result:

\begin{lemma}
\label{lemma1}Let $S_{W}$ be a free symplectic matrix (\ref{free}). We have 
\begin{equation}
\det (S_{W}-I)=\det B\det (B^{-1}A+DB^{-1}-B^{-1}-(B^{T})^{-1}
\label{bofor1}
\end{equation}%
that is, when $S$ is written in the form (\ref{plq}):%
\begin{equation}
\det (S_{W}-I)=\det (L^{-1})\det (P+Q-L-L^{T})\text{.}  \label{bofor2}
\end{equation}
\end{lemma}

\begin{proof}
We begin by noting that since $B$ is invertible we can write $S-I$ as%
\begin{equation*}
\begin{bmatrix}
A-I & B \\ 
C & D-I%
\end{bmatrix}%
=%
\begin{bmatrix}
0 & B \\ 
I & D-I%
\end{bmatrix}%
\begin{bmatrix}
C-(D-I)B^{-1}(A-I) & 0 \\ 
B^{-1}(A-I) & I%
\end{bmatrix}%
\end{equation*}%
hence%
\begin{equation*}
\det (S_{W}-I)=\det B\det (C-(D-I)B^{-1}(A-I))\text{.} 
\end{equation*}%
Since $S$ is symplectic we have $C-DB^{-1}A=-(B^{T})^{-1}$ (use for instance
the fact that $S^{T}JS=SJS^{T}=J$) and hence%
\begin{equation*}
C-(D-I)B^{-1}(A-I))=B^{-1}A+DB^{-1}-B^{-1}-(B^{T})^{-1}\text{;} 
\end{equation*}
the Lemma follows.
\end{proof}

\begin{proposition}
Let $S$ be a free symplectic matrix (\ref{free}) and $\widehat{R}(S)$ the
corresponding Mehlig--Wilkinson operator. We have $\widehat{R}(S)=\widehat{S}%
_{W,m}$ provided that $\nu $ is chosen so that%
\begin{equation}
\nu \equiv m-\limfunc{Inert}(P+Q-L-L^{T})\text{ \ }\func{mod}4
\label{Maslov1}
\end{equation}%
($\limfunc{Inert}(P+Q-L-L^{T})$ the number of $<0$ eigenvalues of the
symmetric matrix $P+Q-L-L^{T}$).
\end{proposition}

\begin{proof}
Recall that we have shown that $\widehat{R}(S)=c_{S}\widehat{S}_{W,m}$ where 
$c_{S}$ is a complex constant with $|c_{S}|=1$. Let us determine that
constant. Let $\delta \in \mathcal{S}^{\prime }(\mathbb{R}^{n})$ be the
Dirac distribution centered at $x=0$; setting%
\begin{equation*}
C=\left( \frac{1}{2\pi }\right) ^{n}\frac{i^{\nu }}{\sqrt{|\det (S_{W}-I)}} 
\end{equation*}%
we have, by definition of $\widehat{R}(S)$, 
\begin{eqnarray*}
\widehat{R}(S)\delta (x) &=&C\int e^{\frac{i}{2}\left\langle
M_{S}z_{0},z_{0}\right\rangle }e^{i(\left\langle p_{0},x\right\rangle -\frac{%
1}{2}\left\langle p_{0},x_{0}\right\rangle )}\delta (x-x_{0})d^{2n}z_{0} \\
&=&C\int e^{\frac{i}{2}\left\langle M_{S}(x,p_{0}),(x,p_{0})\right\rangle
}e^{\frac{i}{2}\left\langle p,x\right\rangle }\delta (x-x_{0})d^{2n}z_{0}
\end{eqnarray*}%
hence, setting $x=0$,%
\begin{equation*}
\widehat{R}(S)\delta (0)=C\int e^{\frac{i}{2}\left\langle
M_{S}(0,p_{0}),(0,p_{0})\right\rangle }\delta (-x_{0})d^{2n}z_{0} 
\end{equation*}%
that is, since $\int \delta (-x_{0})d^{n}x_{0}=1$,%
\begin{equation}
\widehat{R}(S)\delta (0)=\left( \frac{1}{2\pi }\right) ^{n}\frac{i^{\nu }}{%
\sqrt{|\det (S-I)}}\int e^{\frac{i}{2}\left\langle
M_{S}(0,p_{0}),(0,p_{0})\right\rangle }d^{n}p_{0}\text{.}  \label{sdo}
\end{equation}%
Let us calculate the scalar product 
\begin{equation*}
\left\langle M_{S}(0,p_{0}),(0,p_{0})\right\rangle =\sigma
((S-I)^{-1}0,p_{0}),(0,p_{0}))\text{.} 
\end{equation*}%
The relation $(x,p)=(S-I)^{-1}(0,p_{0})$ is equivalent to $%
S(x,p)=(x,p+p_{0}) $ that is to%
\begin{equation*}
p+p_{0}=\partial _{x}W(x,x)\text{ \ and \ }p=-\partial _{x^{\prime }}W(x,x)%
\text{.} 
\end{equation*}%
Using the explicit form (\ref{wplq}) of $W$ together with Lemma \ref{lemma1}
these relations yield%
\begin{equation*}
x=(P+Q-L-L^{T})^{-1}p_{0}\text{ \ ; \ }p=(L-Q)(P+Q-L-L^{T})^{-1}p_{0} 
\end{equation*}%
and hence%
\begin{equation}
\left\langle M_{S}(0,p_{0}),(0,p_{0})\right\rangle =-\left\langle
(P+Q-L-L^{T})^{-1}p_{0},p_{0}\right\rangle \text{.}  \label{bofor3}
\end{equation}%
Applying Fresnel's formula (\ref{fres}) we get%
\begin{equation*}
\left( \frac{1}{2\pi }\right) ^{n}\int e^{\frac{i}{2}\left\langle
M_{S}(0,p_{0}),(0,p_{0})\right\rangle }d^{n}p_{0}=e^{-\frac{i\pi }{4}%
\limfunc{sgn}(P+Q-L-L^{T})}|\det (P+Q-L-L^{T})|^{1/2}\text{;} 
\end{equation*}%
since 
\begin{equation*}
\frac{1}{\sqrt{|\det (S-I)}}=|\det L|^{1/2}|\det (P+Q-L-L^{T})|^{-1/2} 
\end{equation*}%
in view of (\ref{bofor2}) in Lemma \ref{lemma1} we thus have%
\begin{equation*}
\widehat{R}(S)\delta (0)=\left( \frac{1}{2\pi }\right) ^{n}i^{\nu }e^{-\frac{%
i\pi }{4}\limfunc{sgn}(P+Q-L-L^{T})}|\det L|^{1/2}\text{.} 
\end{equation*}%
Now, by definition of $\widehat{S}_{W,m}$ we have%
\begin{equation*}
\widehat{S}_{W,m}\delta (0)=\left( \frac{1}{2\pi }\right) ^{n}i^{m-n/2}|\det
L|^{1/2} 
\end{equation*}%
hence%
\begin{equation*}
i^{\nu }e^{-\frac{i\pi }{4}\limfunc{sgn}(P+Q-L-L^{T})}=i^{m-n/2}\text{.} 
\end{equation*}%
It follows that we have%
\begin{equation*}
\nu -\frac{1}{2}\limfunc{sgn}(P+Q-L-L^{T})\equiv m-\frac{n}{2}\text{ \ }%
\func{mod}4 
\end{equation*}%
which is the same thing as (\ref{Maslov1}) since $P+Q-L-L^{T}$ has rank $n$.
\end{proof}

\subsection{The general case}

Recall that we established in Lemma \ref{lemma1} the equality%
\begin{equation}
\det (S_{W}-I)=\det L^{-1}\det (P+Q-L-L^{T}).  \label{splq}
\end{equation}%
valid for all free matrices $S_{W}\in Sp(n)$. Also recall that every $%
\widehat{S}\in Mp(n)$ can be written (in infinitely many ways) as a product $%
\widehat{S}=\widehat{S}_{W,m}\widehat{S}_{W^{\prime },m^{\prime }}$. We are
going to show that $\widehat{S}_{W,m}$ and $\widehat{S}_{W^{\prime
},m^{\prime }}$ in addition always can be chosen such that $\det (\widehat{S}%
_{W,m}-I)\neq 0$ and $\det (\widehat{S}_{W^{\prime },m^{\prime }}-I)\neq 0$.
For that purpose we need the following straightforward factorization result
(see \cite{AIF}):

\begin{lemma}
Let $W$ be given by (\ref{wplq}); then 
\begin{equation}
\widehat{S}_{W,m}=\widehat{V}_{-P}\widehat{M}_{L,m}\widehat{J}\widehat{V}%
_{-Q}  \label{fac1}
\end{equation}%
where 
\begin{equation*}
\widehat{V}_{-P}f(x)=e^{\frac{i}{2}\left\langle Px,x\right\rangle }f(x)\text{
\ ; \ }\widehat{M}_{L,m}f(x)=i^{m}\sqrt{|\det L|}f(Lx) 
\end{equation*}%
and $\widehat{J}$ is the modified Fourier transform given by%
\begin{equation*}
\widehat{J}f(x)=\left( \frac{1}{2\pi i}\right) ^{n/2}\int e^{-i\left\langle
x,x^{\prime }\right\rangle }f(x^{\prime })d^{n}x^{\prime }\text{.} 
\end{equation*}
\end{lemma}

Let us now state and prove the main result of this section:

\begin{proposition}
Every $\widehat{S}\in Mp(n)$ is the product of two Mehlig--Wilkinson
operators; these operators thus generate $Mp(n)$.
\end{proposition}

\begin{proof}
Let us write $\widehat{S}=\widehat{S}_{W,m}\widehat{S}_{W^{\prime
},m^{\prime }}$ and apply (\ref{fac1}) to each of the factors; this yields
(with obvious notations)%
\begin{equation}
\widehat{S}=\widehat{V}_{-P}\widehat{M}_{L,m}\widehat{J}\widehat{V}%
_{-(P^{\prime }+Q)}\widehat{M}_{L^{\prime },m^{\prime }}\widehat{J}\widehat{V%
}_{-Q^{\prime }}\text{.}  \label{sprod}
\end{equation}%
We claim that $\widehat{S}_{W,m}$ and $\widehat{S}_{W^{\prime },m^{\prime }}$
can be chosen in such a way that $\det (\widehat{S}_{W,m}-I)\neq 0$ and $%
\det (\widehat{S}_{W^{\prime },m^{\prime }}-I)\neq 0$ that is, 
\begin{equation*}
\det (P+Q-L-L^{T})\neq 0\text{ \ and \ }\det (P^{\prime }+Q^{\prime
}-L^{\prime }-L^{\prime T})\neq 0\text{.} 
\end{equation*}
This will prove the assertion in view of (\ref{splq}). We first remark that
the right hand-side of (\ref{sprod}) obviously does not change if we replace 
$P^{\prime }$ by $P^{\prime }+\lambda I$ and $Q$ by $Q-\lambda I$ where $%
\lambda \in \mathbb{R}$. Choose now $\lambda $ such that it is not an
eigenvalue of $P+Q-L-L^{T}$ and $-\lambda $ is not an eigenvalue of $%
P^{\prime }+Q^{\prime }-L^{\prime }-L^{\prime T}$; then 
\begin{equation*}
\det (P+Q-\lambda I-L-L^{T})\neq 0\text{ \ and \ }\det (P^{\prime }+\lambda
I+Q^{\prime }-L-L^{T})\neq 0\text{.} 
\end{equation*}
\end{proof}


\begin{thebibliography}{9}
\bibitem{Folland} \textsc{G. B. Folland}. Harmonic Analysis in Phase space.
Annals of Mathematics studies, Princeton University Press, Princeton, N.J.,
1989.and the references therein.

\bibitem{Bullsci} \textsc{M. de Gosson.} On half-form quantization of
Lagrangian manifolds and quantum mechanics in phase space. Bull. Sci. Math. 
\textbf{121} (1997) 301--322

\bibitem{AIF} \textsc{M. de Gosson.} Maslov Indices on $Mp(n)$.\ Ann. Inst.
Fourier, Grenoble, \textbf{40}(3) (1990) 537--55

\bibitem{IHP} \textsc{M. de Gosson}. On the classical and quantum evolution
of Lagrangian half-forms in phase space. Ann. Inst. H. Poincar\'{e}, \textbf{%
70}(6) (1999) 547--73

\bibitem{Howe} \textsc{R. Howe.} The Oscillator semigroup. Proc. of Symposia
in Pure Mathematics \textbf{48}, Amer. Math. Soc. (1988) 61--132

\bibitem{MW} \textsc{B. Mehlig and M. Wilkinson.} Semiclassical trace
formulae using coherent states. Ann. Phys. \textbf{18}(10), 6--7 (2001)
541-555.
\end{thebibliography}
\end{document}